\documentclass[11pt]{article}
\usepackage{amssymb,amsmath,latexsym, bbm}

\usepackage{epsfig}

\newtheorem{thm}{Theorem}[section]

\newtheorem{cor}[thm]{Corollary}

\newtheorem{lemma}[thm]{Lemma}
\newtheorem{propn}[thm]{Proposition}
\newtheorem{rem}[thm]{Remark}
\newtheorem{example}[thm]{Example}

\def\qed{{\hfill $\Box$ \bigskip}}

\def\R{\mathbb{R}}

\def\Z{\mathbb{Z}}

\def\P{\mathbb{P}}

\def\E{\mathbb{E}}

\def\r{\rho}

\def\EE{\mathcal{E}}
\def\FF{\mathcal{F}}
\def\NN{\mathcal{N}}
\def\JJ{\mathcal{J}}
\def\LL{\mathcal{L}}

\def\1{\mathbbm{1}}

\def\bar{\overline}

\def\wt{\widetilde}
\def\eps{\varepsilon}
\def\tp {\wt \phi}

\def\bee{\begin{equation}}
\def\eee{\end{equation}}

\makeatletter
\@addtoreset{equation}{section}

\makeatother

\begin{document}
\bibliographystyle{plain}

\title{\Large \bf
Symmetric Jump Processes and their Heat Kernel Estimates
}

\author{{\bf Zhen-Qing Chen}
\thanks{Research partially supported by NSF Grant  DMS-06000206.}
}
\date{}
\maketitle

\begin{abstract} We survey the recent development of the
   DeGiorgi-Nash-Moser-Aronson type  theory for a class of symmetric jump processes
   (or equivalently, a class of symmetric integro-differential operators).
   We focus on the
sharp two-sided estimates for the transition density functions (or
heat kernels) of the processes, a priori H\"older estimate and
parabolic Harnack inequalities for their parabolic functions. In
contrast to the second order elliptic differential operator case,
the methods to establish these properties for symmetric
integro-differential operators are mainly probabilistic.
\end{abstract}



\section{Introduction}

Second order elliptic differential operators and diffusion processes
take up, respectively, an central place in the theory of partial
differential equations (PDE) and the theory of probability. There
are close relationships between these two subjects. For a large
class of second order elliptic differential operators $\LL$ on
$\R^n$, there is a diffusion process $X$ on $\R^n$ associated with
it so that $\LL$ is the infinitesimal generator of $X$, and vice
versa.  The connection between $\LL$ and $X$ can also be seen as
follows. The fundamental solution (also called heat kernel)
  for $\LL$ is the transition density function of $X$.
For example, when $\LL=\frac12 \sum_{i, j=1}^n \frac{\partial}{\partial x_i} \left(a_{ij}(x)
 \frac{\partial}{\partial x_j}\right)$, where $(a_{ij}(x))_{1\leq i, j\leq n}$
  is a measurable $n\times n$
 matrix-valued   function on $\R^n$ that is uniformly elliptic and bounded,
there is a symmetric diffusion $X$ having $\LL$ as its $L^2$-infinitesimal generator.
The celebrated  DeGiorgi-Nash-Moser-Aronson theory tells us that
every   bounded parabolic function of $\LL$ (or equivalently, of $X$) is locally H\"older
continuous and the parabolic Harnack inequality holds for non-negative parabolic functions
of $\LL$. Moreover, $\LL$ has a jointly continuous heat kernel $p(t, x, y)$
with respect to the Lebesgue measure on $\R^n$ that enjoys
the following Aronson's estimate:
there are constants $c_k>0$,
 $k=1, \cdots, 4$, so that
\begin{equation}\label{e:1.1}
 c_1 \, p^c(t, c_2|x-y|) \leq p(t, x, y)\leq c_3\, p^c(t, c_4|x-y|) \qquad \hbox{for }
t>0 \hbox{ and } x, y \in \R^n.
 \end{equation}
 Here
 \begin{equation}\label{e:pc}
  p^c(t, r) := t^{-n/2} \exp ( - r^2/t).
 \end{equation}
See \cite{Str} for some history and a survey on this subject, where a mixture of analytic and probabilistic methods is presented.

Recently there are intense interests in studying discontinuous Markov processes,
due to their importance both in theory and in applications.
Many physical and
economic systems should be and in fact have been successfully
modeled by non-Gaussian jump  processes; see for example,
\cite{B, JW, KSZ, STa} and the references therein.
 The infinitesimal generator of a discontinuous Markov process in $\R^n$
is no longer a differential operator but rather a non-local (or,
integro-differential) operator. For instance, the infinitesimal
generator of an isotropically
 symmetric $\alpha$-stable process in $\R^n$ with
$\alpha \in (0, 2)$ is a fractional Laplacian operator $c\,
\Delta^{\alpha /2}:=- c\, (-\Delta)^{\alpha /2}$. During the past
  several years  there are also many interests from
the theory of PDE (such as singular obstacle problems) to study
non-local operators; see, for example,
 \cite{CSS, S} and the references
therein.

 In this paper, we survey   recent development of the
   DeGiorgi-Nash-Moser-Aronson type  theory for
   the following type of non-local
 (integro-differential) operators $\LL$ on $\R^n$:
 \begin{equation}\label{e:1.2}
  \LL  u(x) =\frac12 \sum_{i, j=1}^n \frac{\partial }{\partial x_i} \left(a_{ij}(x)
 \frac{\partial u(x)}{\partial x_j}\right) +
 \lim_{\eps \downarrow 0} \int_{\{y\in \R^n: \, |y-x|>\eps\}}
 (u(y)-u(x)) J(x, y) dy,
 \end{equation}
where either $(a_{ij}(x))_{1\leq i, j\leq n} $ is identically zero or
$(a_{ij}(x))_{1\leq i, j\leq n}$ is a measurable $n\times n$
matrix-valued measurable function on $\R^n$ that is uniformly
elliptic and bounded, and $J$ is a measurable non-negative
symmetric kernel satisfying
certain conditions.
 Associated with such a non-local operator $\LL$
 is an $\R^n$-valued symmetric jump process  $X$ with jumping  kernel $J(x, y)$
 and with possible diffusive components
 when $(a_{ij}(x))_{1\leq i, j\leq n}$ is non-degenerate.
Note that the jumping kernel $J$ determines a L\'evy system for $X$,
which describes the jumps of the process $X$: for any non-negative
measurable function $f$ on $\R_+ \times \R^n\times \R^n$, $t\geq 0 $,
$x\in \R^n$ and stopping time $T$ (with respect to the minimal admissible filtration of $X$),
\begin{equation}\label{e:levy}
\E_x \left[\sum_{s\le T} f(s,X_{s-}, X_s) \right]= \E_x \left[
\int_0^T \left( \int_{\R^n} f(s,X_s, y) J(X_s,y) dy \right) ds \right].
 \end{equation}
Our   focus will be on sharp two-sided heat kernel estimates  for $\LL$ (or,
equivalently, transition density function estimates
for $X$), as well as  parabolic Harnack inequality and
a priori joint H\"older continuity estimate for parabolic functions
of $\LL$. When $(a_{ij}(x))_{1\leq i, j\leq n}\equiv 0$ and $J(x,
y)=c |x-y|^{-n-\alpha}$ for some $\alpha \in (0, 2)$ in
\eqref{e:1.2}, $\LL$ is a fractional Laplacian $c_1
\Delta^{\alpha/2}$ on $\R^n$ and its associated process $X$ is a
rotationally symmetric $\alpha$-stable process on $\R^n$. Unlike the
Brownian motion case, the explicit formula for
 the density function $p(t, x, y)$ of $X$ with respect to the Lebesgue measure
 is only known for a few special $\alpha$, such as $\alpha=1$. However due to
 the scaling property of $X$, one has
  $$ p(t, x, y)=t^{-n/\alpha} \, p(1, t^{-1/\alpha} x,
 t^{-1/\alpha} y)=t^{-n/\alpha}\, f(t^{-1/\alpha}(x-y))
  \quad \hbox{for } t>0 \hbox{ and } x, y\in \R^n,
 $$
 where $f(z)$ is the density function of the symmetric
 $\alpha$-stable random variable $X_1-X_0$
 in $\R^n$. Using Fourier transform, it is not difficult to show
 (see \cite[Theorem 2.1]{BG}) that $f(z)$ is a continuous strictly
 positive function on $\R^n$ depending on $z$ only through $|z|$
 and that
 $f(z)\asymp |z|^{-n-\alpha}$ at infinity.  Consequently
 \begin{equation}\label{e:1.5}
 p(t, x, y)\asymp t^{-n/\alpha} \left(  1 \wedge \frac{t^{1/\alpha}}{|x-y|}
 \right)^{n+\alpha} \qquad \hbox{on } \R_+\times \R^n \times \R^n.
 \end{equation}
In this paper, for two non-negative functions $f$ and $g$, the
notation $f\asymp g$ means that there are positive constants $c_1$
and $c_2$ so that $c_1g( x)\leq f (x)\leq c_2 g(x)$ in the common
domain of definition for $f$ and $g$.  For $a, b\in \R$, $a\wedge
b:=\min \{a, b\}$ and $a\vee b:=\max\{a, b\}$. However such kind of
simple  argument for \eqref{e:1.5} breaks down for the symmetric $\alpha$-stable-like
processes on $\R^n$ when $(a_{ij}(x))_{1\leq i, j\leq n}\equiv 0$
and $J(x, y)=c(x, y) |x-y|^{-n-\alpha}$ for some $\alpha \in (0, 2)$
and a symmetric function $c(x, y)$ that is bounded between two
positive constants in \eqref{e:1.2}, as in this case, $X$ is no longer
a L\'evy process.

Two-sided heat
kernel estimates for jump processes in $\R^n$ have only been studied
recently. In \cite{kol}, Kolokoltsov obtained two-sided heat kernel
estimates for certain stable-like processes in $\R^n$, whose
infinitesimal generators are   a class of pseudo-differential
operators having smooth symbols.
 Bass and Levin
\cite{BL02b} used a completely different approach to obtain similar
estimates for discrete time Markov chain on $\Z^n$ where the
conductance between $x$ and $y$ is comparable to $|x-y|^{-n-\alpha}$
for $\alpha \in (0,2)$. In \cite{CK1}, two-sided heat kernel
estimates and a scale-invariant parabolic Harnack inequality (PHI in abbreviation) for
symmetric $\alpha$-stable-like processes on $d$-sets are obtained.
 Recently in \cite{CK2},
PHI and two-sided heat kernel estimates are
even established for non-local operators of variable order.
Finite range stable-like processes on $\R^n$ are studied in \cite{CKK}.
 This class of processes is  very natural in applications where
jumps only up to a certain size are allowed. The heat kernel
estimates obtained in \cite{CKK} shows finite range stable-like
processes  behave like discontinuous stable-like processes in small
scale and behave like Brownian motion in large scale. Processes
having such properties  may be useful in applications. For example,
in mathematical finance, it has been observed that even though
discontinuous stable processes provide better representations of
financial data than Gaussian processes (cf. \cite{HPR}), financial
data tend to become more Gaussian over a longer time-scale (see
\cite{M} and the references therein). Our heat kernel estimates in
\cite{CKK} show that finite range stable-like processes have this
type of property. Moreover, finite range stable-like processes avoid
large sizes of jumps which can be considered as impossibly huge
changes of financial data in short time.
 See \cite{BBCK} for  some results on parabolic Harnack
inequality and heat kernel estimate for more general non-local
operators of variable order on $\R^n$, whose jumping kernel is
supported on jump size less than or equal to 1. The
DeGiorgi-Nash-Moser-Aronson type  theory is studied very recently in
\cite{CK3} for diffusions with jumps whose infinitesimal generator
is of  type \eqref{e:1.2} with uniformly elliptic and bounded
diffusion matrix $(a_{ij}(x))_{1\leq i, j\leq n}$ and non-degenerate
measurable jumping kernel $J$.

Quite often we need to consider part process $X^D$ of $X$ killed
upon exit an open set $D\subset \R^n$. When $X$ is a Brownian
motion, the infinitesimal generator of $X^D$ is the Dirichlet
Laplacian $\frac12  \Delta_D$. When $X$ is a rotationally symmetric
$\alpha$-stable process in $\R^n$, the infinitesimal generator of
$X^D$ is a Dirichlet fractional Laplacian $c\, \Delta^{\alpha/2}|_D$
that satisfies zero {\it exterior} condition on $D^c$. Though the
transition density function of Brownian motion has been known for
quite a long time, due to the complication near the boundary, a
complete sharp two-sided estimates on the transition density of
killed Brownian motion in bounded $C^{1,1}$ domains $D$
(equivalently, the Dirichlet heat kernel) have been established only
recently in 2002, see \cite{Zq3} and the references therein. Very
recently in \cite{CKS1}, we have obtained sharp two-sided heat
kernel estimates for Dirichlet fractional Laplacian operator in
$C^{1,1}$ open sets, while in \cite{CKS2} and \cite{CKS3}, we
derived sharp two-sided estimates for transition density functions
of censored stable processes  and relativistic $\alpha$-stable
processes in $C^{1,1}$ open sets, respectively.

The rest of the paper is organized as follows. Heat kernel
estimates, PHI and a priori H\"older estimates for stable-like
processes and mixed stable-like processes on $n$-sets in $\R^n$ is
discussed in Sections 2 and 3, respectively. In Section 4, we deal
with finite range stable-like processes on $\R^n$, while results for
diffusions with jumps are surveyed in Section 5. Sections 6 and 7
are devoted to sharp heat kernel estimates for symmetric stable
processes and censored stable processes in $C^{1,1}$-open sets. To
give  a glimpse of our approach to the
   DeGiorgi-Nash-Moser-Aronson type  theory for non-local operators
   using probabilistic means, we give an outline of the main ideas in our investigation
   for the following three classes of processes:
   symmetric stable-like processes on open $n$-sets in $\R^n$ in Section 2,
   diffusions with jumps on $\R^n$ in Section 5
   and symmetric stable processes in open subsets
   of $\R^n$ in   Section 6.
This paper surveys some recent research   that the author is
involved.
  See Bass \cite{B07} for a survey  for related topics on SDEs with jumps, Harnack
inequalities and H\"older continuity of harmonic functions for
non-local operators, and Chen \cite{C} for a survey (prior to 2000)
on potential theory of symmetric stable processes in open sets.

Throughout this paper, $n\geq 1$ is an integer. We denote by $m$ or
$dx$ the $n$-dimensional Lebesgue measure in $\R^n$, and
$C^1_c(\R^n)$ the space of $C^1$-functions on $\R^n$ with compact
support. For a closed subset $F$ of $\R^n$, $C_c(F)$ denotes the
space of continuous functions with compact support in $F$. For a
Markov process $X$ on a state space $E$ and  a subset $K\subset E$, we
let $\sigma_K:= \inf\{ t\geq 0: \, X_t \in K\}$ and $\tau_K:=\inf\{
t\geq 0: \, X_t \notin K\}$ to denote the first entering and exiting
time of $K$ by $X$.

\medskip

{\bf Acknowledgement}. The author thanks Panki Kim and Renming Song
for comments on an earlier version of this paper.

\section{Stable-like processes}

A Borel subset  $F$ in $\R^n$ with $n\ge 1$ is said to be an $n$-set if
 there exist constants $r_0>0$, $C_2> C_1>0$ so that
\begin{equation}\label{eqn:dset}
C_1\, r^n\le m (B(x,r) )\le C_2 \, r^n \qquad \mbox{for all
}~x\in F,~0<r\leq r_0.
\end{equation}
In this section and the next, $B(x,r):=\{y\in F: |x-y|<r\}$ and
$|\cdot|$ is the Euclidean metric in $\R^n$.  Every uniformly
Lipschitz domain in $\R^n$ is an $n$-set, so is its Euclidean
closure. It is easy to check that the classical von Koch snowflake
domain in $\R^2$ is an open $2$-set. An $n$-set can have very rough
boundary since every $n$-set with a subset having zero Lebesgue
measure removed is still an  $n$-set.

For a closed $n$-set $F\subset \R^n$ and $0<\alpha <2$, define
\begin{eqnarray}
\FF  &=& \left\{ u \in L^2(F, m): \, \int_{F\times F}
\frac{(u(x)-u(y))^2}
{|x-y|^{n+\alpha}} m (dx)m (dy) < \infty \right\}
\label{eqn:form1} \\
\EE(u,v)&=& \frac 12 \int_{F\times F} (u(x)-u(y))(v(x)-v(y)) \frac{
c(x, y) } {|x-y|^{n+\alpha}}m (dx)m (dy)  \label{eqn:form2}
\end{eqnarray}
for $u, \, v\in \FF$,
where $c(x, y)$ is a symmetric function on $F\times F$ that is bounded
between two strictly positive constants $C_4>C_3>0$, that is,
\begin{equation}\label{eqn:1.4}
C_3 \leq c(x, y) \leq C_4 \qquad \hbox{for } m \hbox{-a.e. }
   x, y \in F.
\end{equation}
It is easy to check that $(\EE, \FF)$ is a regular Dirichlet form on
$L^2(F, m)$ and therefore there is an associated $m$-symmetric Hunt
process $X$ on $F$ starting from every point in $F$ except for an
exceptional set that has zero capacity. We call such kind of process
a $\alpha$-stable-like process on $F$. Note that when $F=\R^n$ and
$c(x, y)$ is a constant function, then $X$ is nothing but a
rotationally symmetric $\alpha$-stable process on $\R^n$.

\begin{thm}[{\cite[Theorem 1.1]{CK1}}]\label{T:2.1}
 Suppose that $F\subset \R^n$ is a closed
$n$-set and $0<\alpha<2$. Then   $X$ has a H\"older continuous
transition density function $p (t, x, y)$ with respect to $m$. This
in particular implies that $X$ can be modified to start from every
point in $F$ as a Feller process. Moreover, there are constants
$c_2>c_1>0$ that depend only on $n$,  $\alpha$, and the constants
$C_k$, $k=1, \cdots, 4$ in \eqref{eqn:dset} and \eqref{eqn:1.4},
respectively, such that
\begin{equation}\label{e:2.1}
c_1 \min \left\{ t^{-n/\alpha}, \, \frac t{|x-y|^{n+\alpha}}
\right\} \le p(t, x, y)\le c_2\min \left\{ t^{-n/\alpha}, \,  \frac
t{|x-y|^{n+\alpha}} \right\},
\end{equation}
for all $x,y\in F$ and $0<t\leq 1$.
\end{thm}

If $F$ is a global $n$-set in the sense that \eqref{eqn:dset} holds
for every $r>0$, then the heat kernel estimates in \eqref{e:2.1}
holds for every $t>0$.

Note that in \cite[Theorem 1.1]{CK1}, the dependence of $c_1, c_2$
on $(C_1, \cdots, C_4)$ in Theorem \ref{T:2.1} is stated  for every
$\alpha$ except for the case of $0<\alpha =n<2$. The reason is that
in \cite{CK1}, the on-diagonal estimate (Nash's inequality) for the
case of $\alpha = n<2$ was established by using an interpolation
method. This restriction can be removed by an alternative way to
establish Nash's inequality, see \cite[Theorem 3.1]{CK2}.

The detailed heat kernel estimates such as those in \eqref{e:2.1}
are very useful in the study of sample path properties of the
processes. For example, the following is proved in \cite{CK1}.

\begin{thm}[{\cite[Theorem 1.2]{CK1}}]\label{C:2.2}
Under the assumption of Theorem \ref{T:2.1}, for every $x\in F$,
$\P_x$-a.s., the Hausdorff dimension of $X[0, 1]:=\{X_t: \, 0\leq t
\leq 1 \}$ is \, $\alpha \wedge n$.
\end{thm}

In fact, a much stronger result can be derived from Theorem \ref{T:2.1}.
 The following uniform Hausdorff dimensional
result and boundary trace result are established in Remarks 3.10 and
4.4 of \cite{BCR}, respectively.

\begin{thm}\label{T:2.3} Let $D$ be an open $n$-set in $\R^n$ with $n\geq 2$ and
$X$ be an $\alpha$-stable-like process on $\overline D$. Then for
every $x\in \overline D$,
$$ \P_x \left( \dim_H X(E)=  \alpha
\dim_H E \ \hbox{ for all Borel sets } E \subset \R_+ \right) =1
$$
and $\P_x$-a.s.
$$
 \dim_H (X[0, \infty) \cap \partial D)
 = \max\left\{ 1 - \frac{n-\dim_H \partial D}{\alpha},
\ 0 \right\} .
$$
Here for a time set $E\subset \R_+$, $X(E):=\{X_t: t\in E\}$ and
$\dim_H (A)$ is the Hausdorff dimension of a set $A$.
\end{thm}

\bigskip
The approach  to Theorem \ref{T:2.1} in \cite{CK1} is probabilistic
in nature and is motivated by the work of Bass and Levin
\cite{BL02a, BL02b} on stable-like  processes on $\Z^n$ and on
$\R^n$. However there are  new challenges for stable-like processes
on $n$-sets, as paper \cite{BL02a} deals with (possibly
non-symmetric) semimartingale stable-like processes on $\R^n$, when
restricted to the symmetric processes case, requiring
$c(x,y)=f(x,y-x)$ and $f(x,h)$ be an even function in $h$, while
paper \cite{BL02b} is concerned about the transition density
function estimates for discrete time stable-like Markov chains on
$\Z^n$.

By Nash's inequality and \cite[Theorems 3.1 and 3.2]{BBCK}, there is
a properly exceptional set $\NN\subset F$ and a positive symmetric
function $p(t, x, y)$ defined on $(0, \infty) \times (F\setminus
\NN) \times (F\times \NN)$ so that $p(t, x, y)$ is the density
function for $X_t$ under $\P_x$ for every $x\in F\setminus \NN$,
$$ p(t+s, x, y)=\int_F p(s, x, z)p(t,z, y) m(dz) \quad \hbox{for every }
x, y\in F\setminus \NN \hbox{ and } t>0
$$
and
$$ p(t, x, y)\leq c t^{-n/\alpha} \quad \hbox{for every } t>0 \hbox{ and }
x, y\in F\setminus \NN.
$$
Moreover, there is an $\EE$-nest $\{F_k, k\geq 1\}$ of compact sets
so that $\NN = E\setminus \bigcup_{k\geq 1} F_k$ and that for every
$t>0$ and $y\in F\setminus \NN$, $x\mapsto p(t, x, y)$ is continuous
on each $F_k$. The proof of Theorem \ref{T:2.1} given in \cite{CK1}
relies on the following three key propositions. The first
proposition is a tightness result for $X$.

\begin{propn}[{\cite[Proposition 4.1]{CK1}}]\label{P:2.4}
For each $r_0>0$, $A> 0$ and $0<B<1$, there exists $0<\gamma <1$
such that for every $0<r\leq r_0$,
\[
\P_x \left( \tau_{B(x,\, Ar)}<\gamma \, r^\alpha \right)
\le B \qquad \hbox{for every } x\in F\setminus \NN.
\]
Moreover, the constant $\gamma$ can be chosen to depend only on
$(r_0, A, B, n, \alpha)$ and the constants \break $(C_1, C_2, C_3,
C_4)$ in \eqref{eqn:dset}   and \eqref{eqn:1.4} respectively.
\end{propn}

\begin{propn}[{\cite[Proposition 4.2]{CK1}}]\label{P:2.5} \begin{description}
\item{\rm (i)} For each $a >0$, there exists
$c_1>0$ such that for every $x\in F\setminus \NN$,
\begin{equation}\label{eq:uphite}
\P_x \left( \sigma_{B(y, \, ar)} < r^\alpha \right)
\le c_1 \,  \left( \frac{r}{|x-y|}
\right)^{d+\alpha} \quad \hbox{ for every }
r\in (0, \, 2^{1/\alpha}].
\end{equation}
Moreover, the constant $c_1$ above can be chosen to depend only on
$(a, n,  \alpha)$ and on the constants $(C_1, C_2, C_3, C_4)$ in
\eqref{eqn:dset}   and \eqref{eqn:1.4}, respectively.

\item{\rm
(ii)} For each $a , b>0$, there exists
$c_2>0$ such that
\begin{equation}\label{eq:lowhite}
\P_x \left( \sigma_{B(y, \, ar)}<r^\alpha \right)
\ge c_2 \,  \left( \frac{r}{|x-y|}
\right)^{d+\alpha},
\end{equation}
for every  $r\in (0, \, 2^{1/\alpha} ]$ and such that $|x-y|\ge b\,
r$. Moreover, the constant $c_2$ above can be chosen to depend only
on $(a, b,  n,  \alpha)$ and on the constants $(C_1, C_2, C_3, C_4)$
in \eqref{eqn:dset}   and \eqref{eqn:1.4}, respectively.
\end{description}
\end{propn}

The last key proposition is a parabolic Harnack inequality.
For this we need to introduce space-time process
$Z_s:=(V_s, X_s)$, where $V_s=V_0+s$. The filtration generated
by $Z$ satisfying the usual condition
will be denoted as $\{ \widetilde \FF_s; \, s\geq 0\}$.
The law of the space-time
process $s\mapsto Z_s$ starting from $(t, x)$ will be denoted
as $\P^{(t, x)}$.
We say that  a non-negative Borel  measurable function
$q(t,x)$ on $[0, \infty)\times F$ is {\it parabolic}
in a relatively open subset $D$ of $(0, \infty)\times F$ if
for every relatively compact open subset $D_1$ of $D$,
$q(t, x)=\E^{(t,x)} \left[ q (Z_{\tau_{D_1}}) \right]$
for every $(t, x)\in D_1\cap (0, \infty)\times (F\setminus \NN)$,
where $\tau_{D_1}=\inf\{s> 0: \, Z_s\notin D_1\}$.
It is easy to see that
for each $t_0>0$ and $x_0\in F\setminus \NN$, $q(t,x):=p (t_0-t, x,x_0)$
is parabolic on $[0,t_0)\times F$.

\medskip

For each $R_0>0$, we denote $\gamma_{R_0}:=\gamma(R_0, 1/2, 1/2)<1$
the constant in Proposition \ref{P:2.4} corresponding to $r_0=R_0$
and $A=B=1/2$. For $t\le 1$ and $r\le R_0$, we define
\[  Q_{R_0}(t,x,r):=[t,t+\gamma_{R_0} r^{\alpha}]\times (B(x,r)\cap F\setminus \NN).
\]

\begin{propn}[{\cite[Proposition 4.3]{CK1}}]\label{P:2.6}
For every
$R_0>0$, $0<\delta \le \gamma_{R_0}$,
there exists $c>0$
such that for every $z\in F$,
$0<R\le R_0$ and
every non-negative function $q$ on $[0, \infty)\times F$
that is parabolic and bounded on
$[0,3\gamma_{R_0} R^{\alpha}]\times B(z,R)$,
\[\sup_{(t,y)\in Q_{R_0}(\delta  R^{\alpha},z,R/3)}q(t,y)\le c \,
\inf_{y\in B(z,R/3)}q(0,y).\]
Moreover,
the constant $c$ above can be chosen to depend only
on $(R_0, \delta, n,  \alpha)$
and on the constants $(C_1, C_2, C_3, C_4)$
in (\ref{eqn:dset})  and (\ref{eqn:1.4})
respectively.
\end{propn}

Note that the parabolic Harnack inequality implies the elliptic
Harnack inequality.

With the above three propositions, the heat kernel estimates
for $p(t, x, y)$ in Theorem \ref{T:2.1} can be established
for every $t>0$ and $x, y \in F\setminus \NN$.
That $p(t, x, y)$ is jointly continuous and hence the heat kernel
estimates hold for every $t>0$ and $x, y \in F$
comes from the following theorem.

\begin{thm}[{\cite[Theorem 4.14]{CK1}}]\label{T:2.7}
For every $R_0>0$,
there is a constant $c=c(R_0)>0$ such that for
every $0<R\le R_0$ and every bounded parabolic function
$q$ in $Q_{R_0}(0, x_0, \max \{4, \, 4^{1/\alpha} \} R)$,
$$
|q(s, x) -q(t, y)| \leq c \,  \| q \|_{\infty, R} \,
R^{-\beta} \, \left( |t-s|^{1/\alpha} + |x-y|
\right)^\beta
$$
holds for $(s, x), \, (t, y)\in Q_{R_0}(0, x_0, R)$, where $\| q
\|_{\infty, R}:=\sup_{(t,y)\in [0, \, \gamma_{R_0} \max\{4, \,
4^\alpha \} R^\alpha ] \times (F\setminus \NN) } |q(t,y)|$. In
particular, for the transition density function $p (t, x, y)$ of
$X$, there are constants $c>0$ and $\beta>0$ such that for any
$0<t_0<1$, $t, \, s \in [t_0, \,  2]$ and $(x_i, y_i)\in (F\setminus
\NN) \times (F\setminus \NN)$ with $i=1, 2$,
$$
 |p (s, x_1, y_1) -p (t, x_2, y_2)| \leq c \,
t_0^{-(d+\beta)/\alpha}
\left( |t-s|^{1/\alpha} + |x_1-x_2|+|y_1-y_2| \right)^\beta.
$$
Moreover, the constant $c$ above can be chosen to depend only on
$(R_0, t_0, n,   \alpha)$ and on the constants $(C_1, C_2, C_3,
C_4)$ in \eqref{eqn:dset}   and \eqref{eqn:1.4}.
\end{thm}

\section{Mixed stable-like processes}

 In applications, the stochastic model may have more than
one noises. So it is natural to consider mixed stable-like processes
and a mixture of diffusion and jump-type processes.

Let $F$ be a closed global $n$-set in $\R^n$.
 Let $\phi= \phi_1  \psi $ be a strictly increasing continuous functions on $\R_+$,
where   $\psi$ is non-decreasing function on $[ 0, \infty )$ with
$\psi(r)=1$ for $~0<r\leq 1$ that is either the  constant function 1
on $\R_+$ or there are  constants $c_0>0$, $c_2\geq c_1>0$ and
$\gamma_2\geq \gamma_1>0 $ so that
\begin{equation}\label{e:psi1}
  c_1  e^{\gamma_1r} \leq \psi (r)\leq c_2  e^{\gamma_2r} \qquad
 \hbox{ for every }~ 1<r<\infty,
\end{equation}
with
 \begin{equation}\label{e:psi2}
\psi (r+1) \leq c_0 \psi (r) \qquad \hbox{for every } r\geq 1,
\end{equation}
 and  $\phi_1$ is a strictly increasing function on $[0, \infty)$
with $\phi_1 (0)=0$, $\phi_1 (1)=1$ and satisfies the following:
there exist constants $c_2>c_1>0$, $c_3>0$, and $\beta_2\ge
\beta_1>0$ such that
\begin{eqnarray}\label{e:phi1}
 c_1 \Big(\frac Rr\Big)^{\beta_1} &\leq&
\frac{\phi_1 (R)}{\phi_1 (r)}  \ \leq \ c_2 \Big(\frac
Rr\Big)^{\beta_2}
\qquad \hbox{for every } 0<r<R<\infty,\\
\int_0^r\frac {s}{\phi_1 (s)}\, ds &\leq & c_3 \, \frac{ r^2}{\phi_1
(r)} \qquad  \hskip0.8truein \hbox{for every } r>0.\label{e:phi2}
\end{eqnarray}

\begin{rem}\label{R:3.1} {\rm
Note that condition  \eqref{e:phi1} is equivalent to the existence
of constants $ c_4,c_5>1$ and $L_0>1$ such that for every $r>0$,
$$   c_4\phi_1 (r) \leq \phi_1 (L_0r)
\leq c_5 \, \phi_1 (r).
$$
}\end{rem}

 Denote by $d$   the diagonal of $F\times F$ and $J$ be a
symmetric measurable function  on $F\times F \setminus d$ such that
for every $(x, y)\in F\times F\setminus d$,
\begin{equation}\label{jsigm2}
\frac{c_1}{ |x-y|^n \, \phi (c_2|x-y|)}\le J(x, y) \le \frac{c_3} {
|x-y|^n \, \phi (c_4|x-y|)}.
\end{equation}

For $u\in L^2(F, m)$, define $\FF:=\left\{u\in L^2(F; m): \,
\int_{F\times F} (u(x)-u(y))^2 J(x, y) m(dx) m (dy)<\infty\right\}$
and
 \begin{equation}\label{eqn:DE}
   \EE (u, u):=
\int_{F\times F} (u(x)-u(y))(v(x)-v(y))J(x, y) m(dx) m (dy) \quad
\hbox{for } u, v\in \FF.
\end{equation}
For $\beta >0$,
$$  \EE_\beta (u, u):= \EE (u, u) + \beta \int_F u(x)^2 m (dx).
$$

It is not difficult to show that $(\EE, \FF)$ is a regular Dirichlet
form on $L^2(F, m)$ (see \cite[Proposition 2.2 and Remark
4.10(ii)]{CK2}. So there is a symmetric Hunt process $Y$ associated
with it, starting from quasi-every point in $F$. However the next
theorem, which is a special case of \cite[Theorem 1.2]{CK2} (cf.
\cite[Remark 4.4(iv)]{CKK}), says that $X$ can be refined to start
from every point in $F$. Moreover, it has a jointly continuous
transition density function $p(t, x, y)$ with respect to the
Lebesgue measure on $F$.
 The inverse function of the strictly increasing function $t\mapsto
\phi (t)$ is denoted by $\phi^{-1}(t)$.

\medskip

\begin{thm}[Theorem 1.2 of \cite{CK2}]\label{T:3.2}
 Under the above conditions,
there is a conservative Feller process $Y$ associated with $(\EE,
\FF)$ that starts from every point in $F$. Moreover the process $Y$
has a continuous transition density function on
 $(0, \infty )\times F\times F$ with respect to the measure $m$,
 which has   the
following estimates. There are positive constants  $c_1>0$,  $c_2>0$
and $C\geq 1$ such that
\[C^{-1}\left(\frac 1{ \phi^{-1}(t)^n}\wedge \frac{t}{|x-y|^n \, \phi(c_1|x-y|)}\right)
\le p(t,x,y)\le C \left(\frac 1{ \phi^{-1}(t)^n}\wedge
\frac{t}{|x-y|^n \phi(c_2 |x-y|)}\right),\] for every  $t\in (0,1]$
and $x,y\in F$. Moreover, when $\psi\equiv 1$,  the above heat
kernel estimates hold for every $t>0$ and $x,y\in F$.
\end{thm}

 We now give some examples such that Theorem \ref{T:3.2}
applies.

\begin{example}\label{E:3.3}
{\rm  If there is $0<\alpha_1<\alpha_2<2$ and a probability
measure $\nu $
 on $[\alpha_1, \alpha_2]$
such that
$$ \phi (r):= \left( \int_{\alpha_1}^{\alpha_2} r^{-\alpha} \, \nu (d \alpha)\right)^{-1},
$$
then conditions \eqref{e:phi1}-\eqref{e:phi2} are satisfied with
$\psi \equiv 1$. Clearly, $\phi$ is a continuous strictly
increasing function with $\phi (0)=0$ and $\phi (1)=1$. The
condition \eqref{e:phi1} is satisfied with $\gamma \equiv 1$
because
$$ \frac{1}{2^{\alpha_1}}\leq  \frac{\phi(r)}{\phi (2r)} \leq \frac{1}{2^{\alpha_2}}
\qquad \hbox{for any } r>0.
$$
For $r>0$, by Fubini's theorem,
$$ \int_0^r \frac{s}{\phi (s)} ds = \int_0^r \int_{\alpha_1}^{\alpha_2}
r^{1-\alpha} \nu (d\alpha) ds = \int_{\alpha_1}^{\alpha_2}
\frac{1}{2-\alpha} r^{2-\alpha} \nu (d\alpha) \leq \frac1{2-\alpha_2}\,
 \frac{r^2}{\phi (r)}
 $$
and so condition \eqref{e:phi2} is satisfied.
 In this case,
$$ J(x, y) \asymp \int_{\alpha_1}^{\alpha_2}
\frac1{|x-y|^{n+\alpha }} \, \nu (d \alpha ).
$$
 A particular case is
when $\nu$ is a discrete measure. For example,
$\nu$ is a discrete measure concentrate on $\alpha, \beta \in (0, 2)$.
In this case, $J(x, y)= \frac{c_1(x, y)}{|x-y|^{n+\alpha}}+ \frac{c_2(x, y)}{|x-y|^{n+\beta}}$, where $c_i(x, y)$ are two symmetric functions that
are bounded between two positive constants, and
$$\phi (r)\asymp
\min \left\{ r^\alpha, \, r^\beta \right\},  \qquad  \phi^{-1}(r)
\asymp \max \left\{r^{1/\alpha}, \, r^{1/\beta}\right\}.
 $$
 Theorem \ref{T:3.2} gives the precise heat kernel estimates for mixed stable-like
 processes on $F$.
When $F=\R^n$,
 Theorem \ref{T:3.2} in particular gives the heat
kernel estimate for L\'evy processes on $\R^n$ which are linear
combinations of  independent symmetric $\alpha$-stable processes. Of
course, Theorem \ref{T:3.2} holds much more generally, even in the
case of  $F=\R^n$. \qed }
\end{example}

\begin{example}\label{E:3.4}
{\rm Let $Y=\{Y_t, t\geq 0\}$ be the relativistic $\alpha$-stable
processes on $\R^n$ with mass $m_0>0$. That is, $\{Y_t, t\geq 0\}$
is a L\'evy process on $\R^n$ with
\[\E [\exp (i\langle\xi, Y_t-Y_0 \rangle)]=\exp \left(t\left(
m_0^\alpha-(|\xi|^2+m_0^{2})^{\alpha/2} \right) \right).
  \]
 where $\alpha\in (0,2)$. It is shown in \cite{CS} that the corresponding
jumping intensity satisfies
\[J(x,y)\asymp \frac {\Psi(m_0 |x-y|)}{|x-y|^{n+\alpha}},\]
where $\Psi(r)\asymp e^{-r}(1+r^{(n+\alpha-1)/2})$ near $r=\infty$,
and $\Psi(r)=1+\Psi''(0)r^2/2+o(r^4)$ near $r=0$. So conditions
\eqref{e:psi1}-\eqref{e:phi2}) are satisfied with $\gamma_1>0$ for
the jumping intensity kernel for every relativistic $\alpha$-stable
processes on $\R^n$.

When $\alpha=1$, the process is called a relativistic Hamiltonian
process. In this case, the heat kernel can be written as
\[p(t,x,y)=\frac t{(2\pi)^n\sqrt{|x-y|^2+t^2}}\int_{\R^n}
e^{m_0t}e^{-\sqrt{(|x-y|^2+t^2)(|z|^2+m_0^2)}}dz,\]
 see \cite{HS}.
 For simplicity, take $m_0=1$. It can be shown
 that for every $t>0$ and $(x, y)\in \R^n\times
\R^n$,
 \begin{eqnarray*}
&& \frac{c_1 t}{(|x-y|+t)^{n+1}}  \left(  1 \vee
\left(|x-y|+t\right)^{d/2}\right) e^{-c_2 \frac{|x-y|^2}
{\sqrt{|x-y|^2+t^2} }} \\
&\leq & p(t, x, y) \\
 &\leq & \frac{c_3 t}{(|x-y|+t)^{n+1}}
\left(  1 \vee \left(|x-y|+t\right)^{d/2}\right) e^{-c_4
\frac{|x-y|^2} {\sqrt{|x-y|^2+t^2} }}.
\end{eqnarray*}
This in particular implies that for every fixed $t_0>0$, there exist
$c_1,\cdots, c_4>0$ which depend on $t_0$ such that
$$
 c_1\left(t^{-n}\wedge \frac t{|x-y|^{n+1}}\right) e^{-c_2  |x-y| }  \leq  p(t, x, y)
 \leq c_3 \left( t^{-n}\wedge \frac t{|x-y|^{n+1}} \right)
e^{-c_4  |x-y| }
$$
for every $t\in (0, t_0]$ and  $x,y\in \R^n$, which is a special
case of Theorem \ref{T:3.2}. } \qed
\end{example}

 The following
construction of Meyer \cite{Mey75} for jump processes played an important role
in our approach in \cite{CK2}.
  Suppose we have two jump intensity kernels $J(x,y)$ and $ J_0(x,y)$
  on $F\times F$ such that their corresponding pure jump Dirichlet forms
  given in terms of  (\ref{eqn:DE}) with $\FF= \overline{{\cal
  D}(\EE)}^{\EE_1}$
  are regular on $F$. Let $Y=\{Y_t, t\geq 0, \, \P_x, x\in F\setminus \NN\}$ and
$Y^{(0)}=\{Y^{(0)}_t, t\geq 0, \, \P_x, x\in F\setminus \NN_0\}$ be
the processes corresponding to the Dirichlet forms whose L\'evy
densities are $J(x,y)$ and $ J_0(x,y)$, respectively.
  Here $\NN$ and $\NN_0$ are the
properly exceptional sets of $Y$ and $Y^{(0)}$, respectively.
Suppose that
  $J_0(x,y)\le J(x,y)$ and
\[{\cal J}(x):=\int_F(J(x,y)-J_0(x,y)) m (dy)\le c,\]
for all $x\in F$. Let
\begin{equation}\label{qdef}
 J_1(x, y):=J(x, y)-J_0(x, y) \qquad \hbox{and} \qquad
  q(x,y) = \frac{J_1(x, y)}{\JJ(x)}.
\end{equation}
  Then we can construct a process $Y$
corresponding to the jump kernel $J$ from $Y^{(0)}$ as follows. Let
$S_1$ be an exponential random variable of parameter 1 independent
of $Y^{(0)}$, let $C_t=\int _0^t \JJ(Y^{(0)}_s)\, ds$,
and let $U_1$ be
the first time that $C_t$ exceeds $S_1$. We let $Y_s = Y^{(0)}_s$
for $0\le s < U_1$.

At time $U_1$ we introduce a jump from $Y_{U_1-}$ to $Z_1$, where
$Z_1$ is chosen at random according to the distribution
$q(Y_{U_1-},\,y)$. We set $Y_{U_1}=Z_1$, and repeat, using an
independent exponential $S_2$, etc.  Since $\JJ(x)$ is bounded, only
finitely many new jumps are introduced in any bounded time interval.
In \cite{Mey75} it is proved that the resulting process corresponds
to the kernel $J$. See also \cite{INW}. Note that if $\NN_0$ is the
properly exceptional set corresponding to $Y^{(0)}$, then this
construction gives that the properly exceptional set $\NN$ for $Y$
can be chosen to be a subset of $\NN_0$.

Conversely, we can also remove a finite number of jumps from a
process $Y$ to obtain a new process  $Y^{(0)}$. For simplicity,
assume that $J_0(x, y) J_1 (x, y)= 0$.
Suppose one starts with the
process $Y$ (associated with $J$), runs it until the stopping time
$S_1=\inf\{t: J_1(Y_{t-},Y_t)>0\}$, and at that time restarts $Y$ at
the point $Y_{S_1-}$. Suppose one then repeats this procedure over
and over. Meyer \cite{Mey75} proves that the resulting process
$Y^{(0)}$ will correspond to the jump kernel $J_0$. In this case
$\NN_0\subset \NN$.

 Assume that the processes $Y$ and $Y^{(0)}$ have transition
density functions $p(t,x,y)$ and $ p^{(0)}(t,x,y)$, respectively.
Let $\{ \FF_t\}_{t\geq 0}$ be the filtration generated by the
process $Y^{(0)}$.  The following   lemma is shown in \cite[Lemma
2.4]{BBCK} and in \cite[Lemma 3.2]{BGK}.

\begin{lemma}\label{L:3.5}
{\rm (i)} For any $A\in \FF_t$,
\[\P_x \left( \{Y_s = Y^{(0)}_s \hbox{ for all } 0\le s \le t\} \cap A \right) \ge
e^{- t \|{\cal J}\|_\infty }\,   \P_x(A).
\]
{\rm (ii)}  If
$\|J_1\|_\infty<\infty$, then
\[p(t,x,y)\le p^{(0)}(t,x,y)+t\|J_1\|_\infty.\]
\end{lemma}

The use of Lemma \ref{L:3.5} can also significantly simplify the
proofs in \cite{CK1} for results in the last section. The relation
between $Y$ and $Y_0$ can be viewed as the probabilistic counterpart
of the Trotter's semigroup perturbation method. For example, the
proof for Propositions \ref{P:2.5}-\ref{P:2.6} can be simplified by
using Lemma \ref{L:3.5}. See the proof of Propositions 4.9 and 4.11
of \cite{CK2} in this regard.

Comparing with Theorem \ref{T:2.1}, Theorem \ref{T:3.2} says that
the rate function for stable processes of mixed type associated with
\eqref{jsigm2}-\eqref{eqn:DE} is $\phi$. Parabolic Harnack
inequality and a prior H\"older estimate also hold for parabolic
functions of $X$, with this rate function $\phi$. For each $r,t>0$,
we define
\[   Q(t,x,r) :=[t,t+\gamma \phi(r)]\times \left(B(x,r)\cap F\right).
\]

\begin{thm}[{\cite[Theorem 4.12]{CK2}}]\label{T:3.6}
For every $0<\delta \le \gamma$,
there exists $c_1>0$
such that for every $z\in F$,
$R\in (0,1]$ (resp. $R>0$ when $\gamma_1=\gamma_2=0$) and
every non-negative function $h$ on $[0, \infty)\times F$
that is parabolic and bounded on
$[0,\gamma \phi(2R)]\times B(z,2R)$,
\[\sup_{(t,y)\in Q(\delta \phi(R),z,R)}h(t,y)\le c_1 \,
\inf_{y\in B(z,R)}h(0,y).\]
In particular, the following holds for $t\le 1$ (resp. $t>0$
when $\gamma_1=\gamma_2=0$).
\begin{equation}\label{harpt}
\sup_{(s,y)\in
Q((1-\gamma)t,z,\phi^{-1}(t))}
p(s,x,y)\le c \,
\inf_{y\in B(z, \phi^{-1}(t))}p((1+\gamma)t,x,y).
\end{equation}
\end{thm}

\begin{propn}[{\cite[Proposition 4.14]{CK2}}]\label{T:3.7}
For every $R_0\in (0,1]$ (resp. $R_0>0$ when $\gamma_1=\gamma_2=0$),
there are constants $c=c(R_0)>0$ and $\kappa>0$ such that for
every $0<R\le R_0$ and every bounded parabolic function
$h$ in $Q(0, x_0, 2R)$,
$$
|h(s, x) -h(t, y)| \leq c \,  \| h \|_{\infty, F} \,
R^{-\kappa} \, \left( \phi^{-1}(|t-s|) + \r(x,y)
\right)^\kappa
$$
holds for $(s, x), \, (t, y)\in Q(0, x_0, R)$,
where $\|h\|_{\infty, F}:=\sup_{(t,y)\in [0, \, \gamma
\phi(2R) ] \times F } |h(t,y)|$.
In particular, for the transition density function $p (t, x, y)$
of $X$, for any $t_0\in (0,1)$
(resp. any $T>0$ and any $t_0\in (0,T)$ when $\gamma_1=\gamma_2=0$),
there are constants $c=c(t_0)>0$ and $\kappa>0$ such that
for any $t, \, s \in [t_0, \,  1]$ (resp. $t,s\in [t_0,T]$) and $(x_i, y_i)\in
F\times F$ with $i=1, 2$,
$$
 |p (s, x_1, y_1) -p (t, x_2, y_2)| \leq c \,
\frac 1{ \phi^{-1}(t_0)^n\, \phi^{-1}(t_0)^{\kappa}} \left(
\phi^{-1}(|t-s|) + \r(x_1,x_2)+\r(y_1,y_2) \right)^\kappa.
$$
\end{propn}

\section{Finite range stable-like processes}

A finite range $\alpha$-stable-like process $X$ on $\R^n$ is a
symmetric Hunt process on $\R^n$ of purely discontinuous type whose
jumping kernel is  $J(x, y)=\frac{c(x, y)}{|x-y|^{n+\alpha}}
\1_{\{|x-y|\leq \kappa \}}$, where $\alpha \in (0, 2)$, $\kappa >0$
and $c(x, y)$ is a symmetric function on $\R^n\times \R^n$ that is
bounded between two positive constants. The Dirichlet form $(\EE,
\FF)$ associated with $X$ on $L^2(\R^n, m)$ is given by
\begin{eqnarray} \FF  &=& \left\{ u \in L^2(\R^n; m):
\, \int_{\R^n\times \R^n} \frac{(u(x)-u(y))^2} {|x-y|^{n+\alpha}}
 \1_{\{|x-y|\leq \kappa \}} m (dx)m (dy) < \infty \right\}
\label{e:4.1} \\
&=& \left\{ u \in L^2(\R^n; m): \, \int_{\R^n\times \R^n}
\frac{(u(x)-u(y))^2} {|x-y|^{n+\alpha}}
  m (dx)m (dy) < \infty \right\} , \nonumber \\
 \EE(u,v)&=& \frac 12 \int_{F\times F} (u(x)-u(y))(v(x)-v(y))
\frac{ c(x, y) } {|x-y|^{n+\alpha}}  \1_{\{|x-y|\leq \kappa \}}
 m (dx)m (dy)  \label{e:4.2}
\end{eqnarray}
for $u, \, v\in \FF$.
 The $L^2$-infinitesimal generator of  $X$ and $(\EE, \FF)$ is  a non-local
 (integro-differential) operators $\LL$ on $\R^n$ with measurable symmetric
 kernel $J(x, y)=\frac{c(x, y)}{|x-y|^{n+\alpha}}
\1_{\{|x-y|\leq \kappa \}}$:
 $$ \LL u(x) =\lim_{\eps \downarrow 0} \int_{\{y\in \R^n: \, |y-x|>\eps\}}
 (u(y)-u(x)) J(x, y) dy.
 $$

\begin{thm}[{\cite[Proposition 2.1 and Theorems 2.3 and 3.6]{CKK}}]
\label{T:4.1} The finite range stable-like process $X$ has a jointly
continuous transition density function $p(t, x, y)$ and so $X$ can
be refined to start from every point on $\R^n$. Moreover the
following sharp two-sided heat kernel estimates hold.

\begin{description}
\item{{\rm(i)}}
 There is $R_*\in (0, 1)$ so that for every $t\in (0, R_*^\alpha]$ and   $x, y\in \R^n$
with $|x-y|\leq R_*$
\[p(t, x,y)\asymp \left(t^{-n /\alpha}\wedge\frac t{|x-y|^{n+\alpha}}\right);\]

\item{{\rm(ii)}} There exists $C_*\in (0, 1)$ such that for
 $x, y\in \R^n$ with $|x-y|\geq \max\{t/C_*, R_*\}$,
\[p(t, x,y)\asymp \left(\frac t{|x-y|}\right)^{c |x-y|}
=\exp \left(-c |x-y|\log \frac {|x-y|}t \right);\]

\item{{\rm(iii)}} For $t\geq R_*^\alpha$ or
$x, y\in \R^n$ with $|x-y|\in [R_*, t/C_*]$,
$$
  p(t, x,y)\asymp t^{-n /2}\exp\left( -\frac{c|x-y|^2}t \right).
$$
\end{description}
\end{thm}

The following weighted Poincar\'e inequality for non-local operators
together with Lemma \ref{L:3.5} played a crucial role in our proof
of Theorem \ref{T:4.1} in \cite{CKK}. In the remainder of this
paper, $B(x, r)$ denotes the Euclidean ball in $\R^n$ with radius
$r$ centered at $x$.

\begin{thm}[{\cite[Proposition 3.2]{CKK}}]\label{T:4.2}
Suppose that $J(x, y)$ is a symmetric non-negative kernel on $\R^n\times \R^n$ such that
 $J(x, y)=0$ when $|x-y|\geq 1$ and
 $$ \kappa_1 |x-y|^{-n-\alpha} \leq J(x, y) \leq \kappa_2 |x-y|^{-n-\beta}
 \qquad \hbox{when } |x-y|< 1
 $$
 for some constants $\kappa_1, \kappa_2>0$ and $0<\alpha<\beta < 2$.
 Let $\phi (x):= c \left(1-|x|^2 \right)^{12/(2-\beta)} \1_{B(0, 1)} (x)$,
where $c>0$ is the normalizing constant so that $\int_{\R^n} \phi
(x) dx =1$.
 Then there  is a positive constant
$c_1=c_1(n, \alpha, \beta)$
independent of $r>1$, such that for
every $u\in L^1(B(0, 1), \phi dx)$,
\begin{eqnarray*}
  \int_{B(0, 1)} ( u(x)-u_{\phi})^2 \phi (x) dx
 \leq   c_1 \int_{B(0, 1)\times B(0, 1)} (u(x)-u(y))^2 \, r^{n+2}
J(rx, r y)  \sqrt{\phi (x) \phi (y)}\, dxdy  .
\end{eqnarray*}
Here $u_\phi := \int_{B(0, 1)} u(x) \phi (x) dx$.
\end{thm}

\section{Diffusions with jumps}

In this section, we consider symmetric Markov processes on $\R^n$
that have both the diffusion and pure jumping components.
More precisely, consider the following regular
 Dirichlet form  $(\EE, \FF )$
  on $L^2(\R^n; m)$ given by
 \begin{equation}\label{e:DF}
\begin{cases} \displaystyle
\EE (u,v)=\frac12 \int_{\R^n} \nabla u(x)\cdot A(x) \nabla
v(x)dx+\int_{\mathbb R^n}
(u(x)-u(y))(v(x)-v(y))J(x,y)dxdy , \\
\hskip 0.3truein \FF  \, = \, \overline{C^1_c (\R^n)}^{\EE_1} ,
\end{cases}
\end{equation}
where $A(x)=(a_{ij}(x))_{1\leq i, j\leq n}$ is a measurable $n\times
n$ matrix-valued function on $\R^n$ that is uniform elliptic and
bounded in the sense that there exists a constant $c\geq 1$ such
that
\begin{equation}\label{unielli}
c^{-1}\sum_{i=1}^n\xi_i^2 \le \sum_{i,j=1}^n a_{ij}(x)\xi_i\xi_j\le
c\sum_{i=1}^n\xi_i^2 \qquad \hbox{for every  } x,
(\xi_1,\cdots,\xi_d)\in \R^n,
\end{equation}
and $J$ is a symmetric non-negative measurable kernel on $\R^n\times
\R^n$
  such that there are positive constants $\kappa_0>0 $,
 and $\beta \in (0, 2)$ so that
\begin{equation}\label{e:J1}
  J(x, y) \leq \kappa_0
|x-y|^{-n-\beta } \qquad \hbox{for } |x-y|
 \leq \delta_0,
 \end{equation}
and that
\begin{equation}\label{e:J2}
\sup_{x\in \R^n} \int_{\R^n} (|x-y|^2 \wedge 1) J(x, y)\, dy <\infty
.
\end{equation}
 Clearly under condition \eqref{e:J1}, condition \eqref{e:J2} is equivalent to
$$ \sup_{x\in \R^n} \int_{\{y\in \R^n: |y-x|\geq 1\}}  J(x, y) \, dy <\infty .
$$
By the Dirichlet form theory, there is an $\R^n$-valued symmetric Hunt process $X$ associated with $(\EE, \FF)$.
The $L^2$-infinitesimal generator of $X$ ia
a non-local (pseudo-differential) operators $\LL $ on $\R^n$:
 \begin{equation}\label{e:op}
  \LL  u(x) =\frac12 \sum_{i, j=1}^n \frac{\partial}{\partial x_i} \left(a_{ij}(x)
 \frac{\partial u(x)}{\partial x_j}\right) +
 \lim_{\eps \downarrow 0} \int_{\{y\in \R^n: \, |y-x|>\eps\}}
 (u(y)-u(x)) J(x, y) dy,
 \end{equation}

When the jumping kernel $J\equiv 0$ in \eqref{e:op} and
\eqref{e:DF}, $\LL$ is a uniform elliptic operator of divergence
form and $X$ is a symmetric diffusion on $\R^n$. It is well-known
that $X$ has a joint H\"older continuous transition density function
$p(t, x, y)$, which enjoys the  celebrated Aronson's
two-sided heat kernel estimate \eqref{e:1.1}.

 When $A(x)\equiv 0$ in \eqref{e:DF} and $J$ is given by
  \begin{equation}\label{e:J4}
J(x, y) \asymp \frac1{|x-y|^{d}\, \phi (|x-y|)},
\end{equation}
where $\phi$  a strictly increasing continuous function $\phi:
\R_+\to \R_+$ with $\phi (0)=0$, and $\phi(1)=1$ that satisfies the
conditions \eqref{e:phi1}-\eqref{e:phi2} with $\phi$ in place of
$\phi_1$ there, the corresponding process $X$ is a mixed stable-like
process on $\R^n$ appeared in the previous section.  We know from
Theorem \ref{T:3.2} that  there are positive constants $0<c_1<c_2$
so that
$$ c_1 p^j(t,  |x-y|) \leq p(t, x, y)\leq c_2 p^j(t,  |x-y|) \qquad \hbox{for }
t>0, x, y \in \R^n,
 $$
where
\begin{equation}\label{eqn:4}
 p^j(t, r) :=\left( \phi^{-1}(t)^{-n} \wedge
\frac{t}{r^n \phi (r)} \right)
\end{equation}
with $\phi^{-1}$ being the inverse function of $\phi$.

In this section, we consider the case where both $A$ and $J$ are
non-trivial in \eqref{e:op} and \eqref{e:DF}.  Clearly such
non-local operators and  diffusions with jumps take up an important
place both in
 theory and in applications. However there are very limited
 work in literature for this mixture case
 on the topics of this paper until very recently.
 One of the difficulties in obtaining fine properties for such an operator
 $\LL$ and process $X$ is that they exhibit different scales: the diffusion part
 has Brownian scaling $r\mapsto   r^2$ while the pure jump part has
 a different type of scaling.
 Nevertheless, there is a folklore which says that with the presence of
the diffusion part corresponding to $\frac12
\sum_{i, j=1}^n \frac{\partial}{\partial x_i} \left(a_{ij}(x)
 \frac{\partial}{\partial x_j}\right)$, better results can be expected
 under weaker assumptions on the jumping kernel $J$ as the diffusion
 part helps to smooth things out. Our investigation in \cite{CK3} confirms such an
 intuition. In fact we can establish a priori H\"older estimate and
 parabolic Harnack inequality under weaker conditions than \eqref{e:J4}.
 We now present the main results of \cite{CK3}.
Let $W^{1,2}(\R^n)$ denote the Sobolev space of order $(1, 2)$ on
$\R^n$; that is, $W^{1,2}(\R^n):=\{ f\in L^2(\R^n; m): \, \nabla f
\in
 L^2(\R^n; m)\}$.
It is not difficult (see Proposition 1.1 of \cite{CK3}) to show that
under the conditions
\eqref{unielli}-\eqref{e:J2},
the domain of the Dirichlet form of \eqref{e:DF} is characterized by
$$\FF = W^{1,2}(\R^n)
$$
and that (Theorem 2.2 of \cite{CK3}) the corresponding process $X$ has infinite lifetime.
  Let
  $Z=\{Z_t:=(V_0-t, X_t), t\geq 0\}$
 denote the
space-time process of $X$. We say that  a non-negative real valued
Borel  measurable function $h(t,x)$ on $[0, \infty)\times \R^n$ is
{\it parabolic} (or {\it caloric}) on $D=(a,b)\times B(x_0,r)$ if
there is a properly exceptional set $\NN\subset \R^n$ such that for
every relatively compact open subset $D_1$ of $D$,
$$
h(t, x)=\E^{(t,x)} [h (Z_{\tau_{D_1}})]
$$
for every $(t, x)\in D_1\cap ([0,\infty)\times
  (\R^n\setminus \NN))$,
where $\tau_{D_1}=\inf\{s> 0: \, Z_s\notin D_1\}$. We remark that in
Sections 2 and 3 the space-time process is defined to be $(V_0+t,
X_t)$ but this is merely a notational difference. (For reader's
convenience, we keep the notations same as those in the references
 \cite{CK1, CK2, CK3}.)

\begin{thm}[Theorem 1.2 of \cite{CK3}]\label{T:holder}
Assume that
the Dirichlet form $(\EE, \FF )$ given by \eqref{e:DF}
satisfies the conditions \eqref{unielli}-\eqref{e:J2}
and that
for every $0<r<\delta_0$,
\begin{equation}\label{e:J1b}
\inf_{x_0, y_0\in \R^n \atop |x_0-y_0|=r} \,
 \inf_{x\in B(x_0, \, r/16)} \int_{B(y_0, \, r/16)} J(x,z )dz>0.
 \end{equation}
Then for every $R_0\in (0,1]$,
there are constants $c=c(R_0)>0$ and $\kappa>0$ such that for
every $0<R\le R_0$ and every bounded parabolic function
$h$ in $Q(0, x_0, 2R):=(0, 4R^2) \times B(x_0, 2R)$,
$$
|h(s, x) -h(t, y)| \leq c \,  \| h \|_{\infty, R} \,
R^{-\kappa} \, \left(  |t-s|^{1/2} + |x-y|
\right)^\kappa
$$
holds for $(s, x), \, (t, y)\in
 (3R^2, 4R^2)\times B(x_0, R)$,
where $\|h\|_{\infty, R}:=\sup_{(t,y)\in [0, \, 4R^2 ] \times \R^n
\setminus \NN } |h(t,y)|$. In particular, $X$ has a jointly
continuous  transition density function $p (t, x, y)$ with respect
to the Lebesgue measure. Moreover, for every $t_0\in (0, 1)$ there
are constants $c >0$ and $\kappa>0$ such that
  for any $t, \, s \in (t_0, \,  1]$
and $(x_i, y_i)\in \R^n\times \R^n$ with $i=1, 2$,
$$
 |p (s, x_1, y_1) -p (t, x_2, y_2)| \leq c \,
  t_0^{-(n+\kappa)/2}
\left(  |t-s|^{1/2} + |x_1-x_2|+|y_1-y_2| \right)^\kappa.
$$
\end{thm}

In addition to \eqref{unielli}-\eqref{e:J2} and \eqref{e:J1b}, if there is a constant $c>0$
 such that
\begin{equation}\label{e:J3}
 J(x,y) \le \frac{c}{r^n}\int_{B(x,r)}
  J(z, y) dz
 \quad
\hbox { whenever $r\le \frac 12 |x-y|\wedge 1$},\, x,y\in \R^n,
\end{equation}
 then the following parabolic Harnack principle holds for non-negative
parabolic functions of $X$.

\begin{thm}[Theorem 1.3 of \cite{CK3}]\label{T:PHI}
Suppose that the Dirichlet form $(\EE, \FF )$ given by \eqref{e:DF}
satisfies the condition \eqref{unielli}-\eqref{e:J2}, \eqref{e:J1b}
and \eqref{e:J3}. For every $\delta\in (0, 1)$, there exist
constants $c_1=c_1( \delta)$ and $ c_2=c_2( \delta)>0$ such that for
every $z\in \R^n$, $t_0\ge 0$, $0<R\leq c_1$ and every non-negative
function $u$ on $[0, \infty)\times \R^n$ that is parabolic on
$(t_0,t_0+6\delta R^2 )\times B(z,4R)$,
\begin{equation}\label{eqn:4.1}
\sup_{(t_1,y_1)\in Q_-}u(t_1,y_1)
\le c_2 \,
\inf_{(t_2,y_2)\in Q_+}u(t_2,y_2),
\end{equation}
where $Q_-=(t_0+\delta R^2,t_0+2\delta R^2)\times B(x_0,R)$ and
$Q_+=(t_0+3\delta R^2,t_0+ 4\delta R^2)\times B(x_0,R)$.
\end{thm}

We next present a two-sided heat kernel estimate for $X$ when $J(x, y)$
satisfies the condition \eqref{e:J4}. Clearly
\eqref{e:J1}-\eqref{e:J2}, \eqref{e:J1b}
 and \eqref{e:J3} are satisfied when
\eqref{e:J4} holds. Recall that functions $p^c(t, x, y)$ and $p^j(t, x, y)$
are defined by \eqref{e:pc} and \eqref{eqn:4}, respectively.

\begin{thm}[Theorem 1.4 of \cite{CK3}]\label{mainHK}
Suppose that \eqref{unielli} holds and that the jumping kernel $J$
of the Dirichlet form $(\EE, \FF )$ given by \eqref{e:DF} satisfies
the condition \eqref{e:J4}. Denote by $p(t, x, y)$ the continuous
transition density function of  the symmetric Hunt process $X$
associated with the regular Dirichlet form $(\EE, \FF )$ of
\eqref{e:DF} with the jumping kernel $J$ given by \eqref{e:J4}.
There are positive constants $c_i$, $i=1, 2, 3, 4$ such that for
every $t>0$ and $x, y \in \R^n$,
\begin{eqnarray}
&& c_1\,    \left( t^{-n/2}\wedge\phi^{-1}(t)^{-n} \right)   \wedge
\left( p^c(t, c_2|x- y|)+p^j(t, |x- y|) \right)
\nonumber\\
&\leq & p(t, x, y) \leq c_3\, \left( t^{-n/2}\wedge\phi^{-1}(t)^{-n}
\right) \wedge \left( p^c(t, c_4 |x- y|)+p^j(t, |x- y|) \right).
\label{eq:HKj+dul}
\end{eqnarray}
Here $p^c$ and $p^j$ are the functions given by \eqref{e:pc} and
\eqref{eqn:4}, respectively.
\end{thm}

When $A(x)\equiv I_{n\times n}$, the $n\times n$ identity matrix,
and $J(x, y)= c |x-y|^{-n-\alpha}$ for some $\alpha \in (0, 2)$ in
\eqref{e:DF}, that is, when $X$ is the independent sum of a Brownian
motion $W$ on $\R^n$ and an isotropically symmetric $\alpha$-stable
process $Y$ on $\R^n$, the transition density function $p(t, x, y)$
can be expressed as the convolution of the transition density
functions of $W$ and $Y$,
 whose two-sided estimates are known.
In \cite{SV07},  heat kernel estimates for this L\'evy process $X$
are carried out by computing the convolution and the estimates  are given
in a form that depends on which region the point $(t, x, y)$ falls into.
Subsequently, the parabolic Harnack inequality \eqref{eqn:4.1}
 for such a L\'evy process
$X$ is derived in \cite{SV07} by using the two-sided heat kernel
estimate. Clearly such an approach is not applicable in our setting
even when $\phi (r)=r^\alpha$,
  since in our case, the diffusion and jumping part of $X$ are typically
not independent. The two-sided estimate in this simple form of
(\ref{eq:HKj+dul})
is a new observation of \cite{CK3} even in the independent sum of a
Brownian motion and an isotropically symmetric $\alpha$-stable process case
considered in \cite{SV07}.

The approach in \cite{CK3} employs methods from both probability
theory and analysis, but it is mainly probabilistic. It   uses some
ideas previously developed in \cite{BBCK, BGK, CK1, CK2, CKK}. To
get a priori H\"older estimates for parabolic functions of $X$, we
establish the following three key ingredients.
\begin{description}
\item{(i)} Exit time upper bound estimate:
$$\E_x [ \tau_{B(x_0, r)}] \leq c_1r^2 \qquad \hbox{for } x\in B(x_0, r),
 $$
 where $\tau_{B(x_0, r)}:=\inf\{t>0: X_t\notin B(x_0, r)\}$
  is the first exit time from $B(x_0, r)$ by $X$.

\item{(ii)} Hitting probability estimate:
$$
\P_x \left(X_{\tau_{B(x, r)}} \notin B(x, s) \right)
\leq  \frac{c_2r^2}{(s\wedge 1)^2} \qquad \hbox{for every } r\in (0, 1]
\hbox{ and } s\geq 2r.
$$

\item{(iii)} Hitting probability estimate for space-time process $Z_t=(V_0-t, X_t)$:
 for every $x\in \R^n$, $r\in (0, 1]$ and any compact subset $A\subset Q(x, r):=(0, r^2)\times B(x, r)$,
 $$
  \P^{(r^2, x)}
 (\sigma_A <\tau_r) \geq c_3
 \frac{m_{n+1}(A)}{r^{n+2}},
 $$
 where by slightly abusing the notation,
 $\sigma_A:=\{t>0: Z_t\in A\}$ is the first hitting time of $A$, $\tau_r:=\inf\{t>0: Z_t\notin Q(x, r)\}$
  is the first exit time from $Q(x, r)$ by $Z$
  and $m_{n+1}$ is the
 Lebesgue measure on $\R^{n+1}$.
\end{description}

Here we use the following notations.
The probability law of the process $X$ starting from $x$ is denoted
as $\P_x$ and the mathematical expectation under it is denoted as
$\E_x$, while probability law of the space-time process $Z=(V, X)$
starting from $(t, x)$, i.e. $(V_0,X_0)=(t,x)$, is denoted as
$\P^{(t, x)}$ and the mathematical expectation under it is denoted
as $\E^{(t, x)}$. To establish parabolic Harnack inequality, we need
in addition the following.

\begin{description}\item{(iv)} Short time near-diagonal heat kernel estimate:
  for every $t_0>0$, there is $c_4=c_4(t_0)>0$
such that for every $x_0\in \R^n$ and $t\in (0, t_0]$,
$$ p^{B(x_0, \sqrt{t})} (t, x, y) \geq c_4 t^{-n/2} \qquad \hbox{for }
x, y\in B(x_0, \sqrt{t}/2).
$$
Here $p^{B(x_0, \sqrt{t})}$ is the transition density function for
the part process $X^{B(x_0, \sqrt{t})}$ of $X$ killed upon leaving
the ball $B(x_0, \sqrt{t})$.

\item{(v)}   Let $R\leq 1$ and $\delta <1$.
$Q_1=[t_0+2\delta R^2/3, \, t_0+ 5\delta R^2]\times B(x_0, 3R/2)$,
$Q_2=[t_0+ \delta R^2/3, \, t_0+11\delta R^2/2]\times B(x_0, 2R)$
and define $Q_-$ and $Q_+$ as in Theorem \ref{T:PHI}. Let $h:
[0,\infty)\times \R^n\to \R_+$ be bounded and supported in
$[0,\infty)\times B(x_0,3R)^c$. Then there exists
$c_5=c_5(\delta)>0$ such that
\[\E^{(t_1,y_1)}[h(  Z_{\tau_{Q_1}})]\le c_5\E^{(t_2,y_2)}
[h(  Z_{\tau_{Q_2}})] \qquad \mbox{for }  (t_1,y_1)\in Q_- \mbox{
and }  (t_2,y_2)\in Q_+.\]
\end{description}

The proof of (iv) uses ideas from \cite{BBCK}, where a similar
inequality is established for finite range pure jump process.
However, some difficulties arise due to the presence of the
diffusion part.

The upper bound heat kernel estimate in Theorem \ref{mainHK} is
established by using method of scaling, by Meyer's construction of
the process $X$ based on finite range process $X^{(\lambda)}$, where
the jumping kernel $J$ is replaced by $J(x, y)\1_{\{|x-y|\leq
\lambda\}}$, and by Davies' method from \cite{CKS} to derive an
upper bound estimate for the transition density function of
$X^{(\lambda)}$ through carefully chosen   testing functions. Here
we need to select the value of $\lambda$ in a very careful way that
depends on the values of $t$ and $|x-y|$.

To get the lower bound heat kernel estimate in Theorem \ref{mainHK},
we need a full scale parabolic Harnack principle that extends
Theorem \ref{T:PHI} to all $R>0$ with the scale function $\tp
(R):=R^2 \wedge \phi (R)$ in place of $R\mapsto R^2$ there.  To
establish such a full scale parabolic Harnack principle, we show the
following.

 \begin{description}
 \item{(iii')} Strengthened version of (iii):
  for every $x\in \R^n$, $r>0$ and any compact subset $A\subset
  Q(0,x,r) :=[0,\gamma_0 \tp(r)]\times B(x,r)$,
 $$    \P^{(\gamma_0 \tp (r), x)}
 (\sigma_A <\tau_r) \geq c_3
 \frac{m_{n+1}(A)}{r^n \tp (r)}.
 $$
 Here $\gamma_0$ denotes the constant $\gamma (1/2, 1/2)$ in Proposition
 6.2 of \cite{CK3}.

 \item{(vi)}
 For
 every $\delta \in (0, \gamma_0]$,
 there is a constant $c_6=c_6(\gamma)$
 so that for every $0<R\leq 1$, $r\in (0, R/4]$
 and $(t, x)\in Q(0, z, R/3)$ with
 $0<t\leq \gamma_0 \tp (R/3)-\delta \tp (r)$,
 $$
  \P^{(\gamma_0\tp (R/3), z)}
  (\sigma_{U(t, x, r)}<\tau_{Q(0, z, R)})
 \geq c_6 \frac{r^n \tp (r)}{R^n \tp (R)},
 $$
 where $U(t, x, r):=\{t\}\times B(x, r)$.
\end{description}
With the full scale parabolic Harnack inequality, the lower bound
heat kernel estimate can then be derived once the following estimate
is obtained.

\begin{description} \item{(vii)}
Tightness result: there are constants
$c_7\geq 2$ and  $c_8>0$ such that for every $t>0$ and $x, y\in
\R^n$ with $|x-y|\geq c_7 \tp (t)$,
$$
\P_x \left(X_t\in B \big(y, c_7\tp^{-1}(t)\big) \right) \geq c_8
\frac{t (\tp^{-1}(t))^n}{|x-y|^n \tp (|x-y|)}.
$$
\end{description}

\section{Dirichlet heat kernel estimates for symmetric stable
processes}

Many times one encounters part process $X^D$ of $X$ killed upon exiting
a open set $D$. The infinitesimal generator $\LL^D$ of $X^D$ is
 the infinitesimal generator $\LL$ of $X$ satisfying  Dirichlet boundary or
 zero exterior condition. It is a fundamental problem both in analysis and in
 probability theory to study precise estimate
 for
 the transition density function
  of $X^D$ (or equivalently, the Dirichlet heat kernel of $\LL^D$).
 However due to the complication near the boundary, two-sided estimates on the
transition density of killed Brownian motion in  bounded $C^{1,1}$
domains $D$ (equivalently, the Dirichlet heat kernel) have been
established only recently in 2002; see \cite{Zq3} and the references
therein.
 In this section, we survey the recent result from \cite{CKS1} on sharp two-sided
 estimates on the transition density function $p_D(t, x, y)$ of part process $X^D$ of
a rotationally symmetric $\alpha$-stable process killed upon leaving a $C^{1,1}$ open
 set $D$. The infinitesimal generator of $X^D$ is the fractional
  Laplacian $c\, \Delta^{\alpha/2}|_D$ satisfying
 zero exterior condition on $D^c$.

 Recall that an open set $D$
in $\R^n$ (when $n\geq  2$) is said to be a $C^{1,1}$ open set if
there exist a localization radius $R_0>0$ and a constant
$\Lambda_0>0$ such that for every $z\in\partial D$, there is a
$C^{1,1}$-function $\phi=\phi_z: \R^{n-1}\to \R$ satisfying $\phi
(0)= \nabla\phi (0)=0$, $\| \nabla \phi  \|_\infty \leq \Lambda_0$,
$| \nabla \phi (x)-\nabla \phi (z)| \leq \Lambda_0 |x-z|$, and an
orthonormal coordinate system $CS_z$: $y=(y_1, \cdots, y_{n-1},
y_n):=(\wt y, \, y_n)$ with its origin at $z$ such that
$$
B(z,R_0)\cap D= \left\{ y\in B(0, R_0): \, y_n>\phi (\wt y)
\right\},
$$
where the ball $B(0, R_0)$ on the right hand side is in the
coordinate system $CS_z$.
 The pair $(R_0, \Lambda_0)$ is called the characteristics of the
$C^{1,1}$ open set $D$. We remark that in some literatures, the
$C^{1,1}$ open set defined above is called a {\it uniform}
 $C^{1,1}$ open set as $(R_0, \Lambda_0)$ is
universal for every $z\in \partial D$. For $x\in \R^n$, let
$\delta_{ D}(x)$ denote the Euclidean distance between $x$ and
$D^c$. By a $C^{1,1}$ open set in $\R$ we mean an open set which can
be written as the union of disjoint intervals so that the minimum of
the lengths of all these intervals is positive and the minimum of
the distances between these intervals is positive. Note that a
$C^{1,1}$ open set can be unbounded and disconnected.

\begin{thm}[{\cite[Theorem 1.1]{CKS1}}]\label{T:6.1}
 Let $D$ be a $C^{1,1}$ open
subset of $\R^n$ with $n\geq 1$ and $\delta_D(x)$ the Euclidean
distance between $x$ and $D^c$.
\begin{description}
\item{\rm (i)} For every $T>0$, on $(0, T]\times D\times D$,
$$
p_D(t, x, y)\asymp t^{-n/\alpha} \left(1\wedge
\frac{t^{1/\alpha}}{|x-y|}\right)^{n+\alpha} \left( 1\wedge
\frac{\delta_D(x)}{t^{1/\alpha}} \right)^{\alpha/2} \left( 1\wedge
\frac{\delta_D(y)}{t^{1/\alpha}}\right)^{\alpha/2}.
$$
\item{\rm (ii)} Suppose in addition that $D$ is bounded.
For every $T>0$, there are positive constants $c_1<c_2$ so that on
$[T, \infty)\times D\times D$,
$$
c_1\, e^{-\lambda_1 t}\, \delta_D (x)^{\alpha/2}\, \delta_D
(y)^{\alpha/2} \,\leq\,  p_D(t, x, y) \,\leq\, c_2\, e^{-\lambda_1
t}\, \delta_D (x)^{\alpha/2} \,\delta_D (y)^{\alpha/2} ,
$$
where $\lambda_1>0$ is the smallest eigenvalue of the Dirichlet
fractional Laplacian $(-\Delta)^{\alpha/2}|_D$.
\end{description}
\end{thm}

By integrating the two-sided heat kernel estimates in Theorem
\ref{T:6.1} with respect to $t$, one can easily recover the
following estimate of
 the Green function $  G_D(x, y)=\int_0^\infty
p_D(t, x, y)dt$, initially obtained independently in \cite{CS1} and
\cite{Ku1} when $n\geq 2$.

\begin{cor}[{\cite[Corollary 1.2]{CKS1}}]\label{C:6.2}
 Let $D$ be a bounded $C^{1,1}$-open set in
$\R^n$ with $n \geq 1$. Then on $D\times D$,
$$
G_D(x, y)\, \asymp\, \begin{cases} \frac{1} {|x-y|^{n-\alpha}}
\left(1\wedge \frac{  \delta_D(x)^{\alpha/2}
\delta_D(y)^{\alpha/2}}{ |x-y|^{\alpha}}
\right)  \qquad &\hbox{when }n >\alpha ,  \\
\log \left( 1+ \frac{  \delta_D(x)^{\alpha/2} \delta_D
(y)^{\alpha/2}}{ |x-y|^{\alpha}
}\right)  &\hbox{when } n=1=\alpha , \\
\big( \delta_D(x)  \delta_D (y)\big)^{(\alpha-1)/2} \wedge \frac{
\delta_D(x)^{\alpha/2} \delta_D (y)^{\alpha/2}}{ |x-y|} &\hbox{when
} n=1<\alpha .
\end{cases}
$$
\end{cor}

\bigskip

 Theorem \ref{T:6.1}(i) is established in \cite{CKS1} through Theorems
\ref{t:ub} and \ref{t:lb}, which give the upper bound and lower
bound estimates, respectively. Theorem \ref{T:6.1}(ii) is an easy
consequence of the intrinsic ultracontractivity of the symmetric
$\alpha$-stable process in a bounded $C^{1,1}$ open set.  In fact,
the upper bound estimates in both Theorem \ref{T:6.1} and Corollary
\ref{C:6.2} hold for any domain $D$ with (a weak version of) the
{\it uniform exterior ball condition} in place of the $C^{1,1}$
condition, while the lower bound estimates in  both Theorem
\ref{T:6.1} and Corollary \ref{C:6.2} hold for any domain $D$ with
the {\it uniform interior ball condition} in place of the $C^{1,1}$
condition.

  We say that  $D$ is an open set
satisfying (a weak version of) the uniform exterior ball condition
with radius $r_0>0$ if for every $z\in \partial D$ and $r\in (0,
r_0)$, there is a ball $B^z$ of radius $r$ such that $B^z\subset
\R^n \setminus \overline D$ and $\partial B^z \cap \partial
D=\{z\}$.

\begin{thm}[{\cite[Theorem 2.4]{CKS1}}]\label{t:ub}
Let $D$ be an open set in $\R^n$ that satisfies
 the uniform exterior ball condition with radius
$r_0>0$. For every $T>0$, there exists a positive constant $c=c(T,
r_0, \alpha)$ such that for $t\in (0, T]$ and $x, y\in D$,
\bee\label{e:1}
   p_D(t, x, y) \leq c \,
  t^{-n/\alpha} \left(1\wedge
\frac{t^{1/\alpha}}{|x-y|}\right)^{n+\alpha} \left( 1\wedge
\frac{\delta_D(x)}{t^{1/\alpha}} \right)^{\alpha/2} \left( 1\wedge
\frac{\delta_D(y)}{t^{1/\alpha}}\right)^{\alpha/2}. \eee
\end{thm}

An  open set $D$ is said to satisfy
 the uniform interior ball condition with radius $r_0>0$ in
the following sense: For every $x\in D$ with $\delta_D (x)<r_0$,
there is
 $z_x\in \partial D$
 so that $|x-z_x|=\delta_D(x)$ and $B(x_0,
r_0)\subset D$ for $x_0:= z_x+r_0 (x-z_x)/|x-z_x|$ .  It is
well-known that any (uniform) $C^{1, 1}$ open set $D$ satisfies both
the {\it uniform interior ball condition} and the {\it uniform
exterior ball condition}.

\begin{thm}[{\cite[Theorem 3.1]{CKS1}}]\label{t:lb}
 Assume that $D$ is an open set in $\R^n$ satisfying the uniform
interior ball condition. Then for every $T>0$ there exists a
positive constant $c=c(r_0,\alpha, T)$ such that for all $(t, x,
y)\in (0, T]\times D\times D$,
$$
p_D(t, x, y) \ge c \,
  t^{-n/\alpha} \left(1\wedge
\frac{t^{1/\alpha}}{|x-y|}\right)^{n+\alpha} \left( 1\wedge
\frac{\delta_D(x)}{t^{1/\alpha}} \right)^{\alpha/2} \left( 1\wedge
\frac{\delta_D(y)}{t^{1/\alpha}}\right)^{\alpha/2}.
$$
\end{thm}

\medskip

There are significant differences between obtaining two-sided
Dirichlet heat kernel estimates for the Laplacian and the fractional
Laplacian, as the latter is a non-local operator.
Our approach in \cite{CKS1} is mainly probabilistic. It uses only the following five
ingredients:

\begin{description}
\item{(i)}
the upper bound heat kernel estimate  for the
rotationally symmetric $\alpha$-stable process $X$ in $\R^n$ and the
stable-scaling property of $X$;

\item{(ii)}
the L\'evy system of $X$ that describes how the process  jumps;

\item{(iii)} the mean exit time estimates from annuli and from balls;

\item{(iv)}
the boundary Harnack inequality of $X$ in annuli (when $n\geq 2$)
and in intervals (when $n=1$), and the parabolic Harnack inequality
of $X$;

\item{(v)}
the intrinsic ultracontractivity of $X$ in bounded open sets.
\end{description}

\medskip

The upper bound heat kernel estimate of $X$ on $\R^n$   gives an upper
bound for $p_D(t, x, y)$, while the L\'evy system is the basic tool
used throughout our argument as the symmetric stable process moves
by ``pure jumping". To get the boundary decay rate of $p_D(t, x,
y)$, we use the boundary Harnack inequality and the domain
 monotonicity of the killed stable process $X^D$ in $D$ by comparing
it with certain truncated exterior balls (i.e. annulus) as well as
interior balls. The mean exit time estimate   for an annulus is applied
with the help of the
boundary Harnack inequality   to get the boundary decay
rate in the upper bound heat kernel estimates.
 The two-sided estimates in the ball $B=B(0,
1)$: $\E_x [\tau_B ]\asymp \delta_B(x)^{\alpha/2}$ is used to get
the two-sided estimate   on the first eigenfunction in
balls. The latter  is then used
to get the boundary decay rate for the lower bound estimate in
$p_D(t,x, y)$. The parabolic Harnack inequality allows us to get
pointwise lower bound on $p_D(t, x, y)$ from the integral of
$w\mapsto p_D(t/2, x, w)$ over some suitable region. When $X^D$ is intrinsic ultracontractive,   $p_D(t, x, y)$ is comparable to $c_t
\phi_D(x)\phi_D(y)$ for some $c_t>0$ and a good control is known for
$c_t$ when $t$ is above a certain large $t_0$,
 where $\phi_D$ is the positive first eigenfunction of
$(-\Delta)^{\alpha/2}|_D$, the infinitesimal generator of $X^D$.

Note that the large time heat kernel estimate in Theorem \ref{T:6.1}(ii) requires $D$ to be bounded.
See \cite{CT} for recent results on large time sharp heat kernel estimates for symmetric stable
processes in certain unbounded $C^{1,1}$ open sets.

The approach developed in \cite{CKS1} is quite general in principle
and can be adapted to study heat kernel estimates for other types of
jump processes in open subsets and their perturbations, such as
censored stable processes to be discussed in next section.

\section{Dirichlet heat kernel estimates for censored stable
processes}

Censored $\alpha$-stable  processes in an open subset of $\R^n$
were introduced and studied by Bogdan, Burdzy and Chen in
\cite{BBC}. Fix an open set $D$ in $\R^n$ with $n\geq 1$. Define a
bilinear form $\EE  $ on $C_c^\infty(D)$ by
\bee\label{e:csdf}
 \EE (u, v):= \frac1{2}  \int_D \int_D (u(x)-u(y))(v(x)-v(y)) \frac{c}
 {|x-y|^{n+\alpha}} dxdy , \quad u, v \in C_c^\infty(D),
 \eee
where  $c>0$ is a constant.
Using Fatou's lemma, it is easy to check that the bilinear form $ (
\EE   , C^\infty_c (D))$ is closable in $L^2 (D, dx)$. Let $\FF$ be
the closure of  $C^\infty_c(D)$  under the Hilbert inner product
$\EE  _1:=\EE  +(\,\cdot\, ,\, \cdot\, )_{L^2(D, dx)}.$ As is noted
in \cite{BBC}, $(\EE  , \FF )$ is Markovian and hence a regular
symmetric Dirichlet form on $L^2(D, dx)$, and therefore there is an
associated symmetric Hunt process $X=\{X_t, t\ge 0, \P_x, x\in D\}$
taking values in $D$. The process
$X$ is  called a censored (or resurrected)
$\alpha$-stable process in $D$.

Let $Y$ be a rotationally symmetric $\alpha$-stable process in $\R^n$ with
jumping kernel $c|x-y|^{-n-\alpha}$.
For any open subset $D$ of $\R^n$, we  use $Y^D$ to denote the
subprocess of $Y$ killed upon exiting from $D$. The following result
gives two other ways of constructing a censored $\alpha$-stable
process.

\medskip

\begin{thm}[{\cite[Theorem 2.1 and Remark 2.4]{BBC}}]\label{T:7.1}
  The following processes
have the same distribution:
\begin{description}
\item{(i)}
the symmetric Hunt process $X$ associated with the regular symmetric
Dirichlet form $(\EE , \FF)$ on  $L^2 (D, dx)$;

\item{(ii)}
the strong Markov process $X$ obtained from the killed symmetric
$\alpha$-stable-like process $Y^D$ in $D$ through the
Ikeda--Nagasawa--Watanabe piecing together procedure;

\item{(iii)}
the process $X$ obtained from $Y^D$ through the Feynman-Kac
transform $e^{\int_0^t \kappa_D(Y^D_s) ds}$ with
$$
\kappa_D(x):=\int_{D^c}\frac{c }{|x-y|^{n+\alpha}}dy.
$$
\end{description}
\end{thm}
\medskip

The Ikeda--Nagasawa--Watanabe piecing together procedure mentioned
in (ii) goes as  follows. Let $X_t(\omega) = Y^D_t(\omega)$ for $t <
\tau_D (\omega)$.  If $Y^D_{\tau_D -}(\omega) \notin D$, set
$X_t(\omega) = \partial$ for $t \ge \tau_D (\omega)$.  If
$Y^D_{\tau_D -}(\omega) \in D$, let $X_{\tau_D}(\omega) =
Y^D_{\tau_D -}(\omega)$ and glue an independent copy of $Y^D$
starting from $Y^D_{\tau_D -}(\omega)$ to $X_{\tau_D}(\omega)$.
Iterating this procedure countably many times, we obtain a process
on $D$ which is a version of the strong Markov process $X$; the
procedure works for every starting point in $D$.

For any open $n$-set $D$ in $\R^n$, define
$$
\FF^{\rm ref}:= \left\{ u\in
L^2(D):\int_{D}\int_{D}\frac{(u(x)-u(y))^2} {|x-y|^{n+\alpha}}\,
dxdy<\infty \right\}
$$
and
$$
\EE ^{\rm ref}(u, v):=\frac12\int_{D}\int_{D}(u(x)-u(y))(v(x)-v(y))
\frac{c } {|x-y|^{n+\alpha}}\, dxdy, \quad u,v\in \FF^{\rm ref}.
$$
 As we see from Section 2, the bilinear form $(\EE
^{\rm ref}, \FF^{\rm ref})$ is a regular symmetric Dirichlet form on
$L^2(\overline D, dx)$. The process $\overline X$ on $\overline D$
associated with $(\EE ^{\rm ref}, \FF^{\rm ref})$ is called in \cite{BBC} a
reflected $\alpha$-stable  process on $\overline D$. By Theorem \ref{T:2.1},
 $\overline X$   has a H\"older
continuous transition density  function $\bar p(t, x, y)$ on $(0,
\infty) \times \overline D \times \overline D$ and for every
$T_0>0$, there are positive constants $c_1, c_2$ so that for $ t\in
(0, T_0]$ and $ x, y \in \overline D$,
 \bee\label{e:rhk}
 c_1 \, t^{-n/\alpha} \left(1\wedge
\frac{t^{1/\alpha}}{|x-y|}\right)^{n+\alpha} \leq \bar p(t, x, y)
\leq c_2 \, t^{-n/\alpha} \left(1\wedge
\frac{t^{1/\alpha}}{|x-y|}\right)^{n+\alpha} .
 \eee
This in particular implies that $\overline X$ can start from every
point in $\overline D$. When $D$ is an open $n$-set in $\R^n$, the
censored $\alpha$-stable-like process $X$ can be realized as a
subprocess of $\overline X$ killed upon leaving $D$ (see
\cite[Remark 2.1]{BBC}). It is proved in \cite{BBC} that when $D$ is
a global Lipschitz domain and $\alpha \in (0, 1]$, then $X$ and
$\overline X$ are the same and so $X$ has a sharp two-sided heat
kernel estimate \eqref{e:rhk} in this case. Hence in the following
we will concentrate on the case of $\alpha \in (1, 2)$. The next
theorem gives a sharp two-sided heat kernel estimate for the
transition density function $p_D(t, x, y)$ of censored
$\alpha$-stable process in an $C^{1,1}$ open set with $\alpha \in
(1, 2)$.

\begin{thm}[{\cite[Theorem 1.1]{CKS2}}]\label{T:7.2}
Suppose that $n\geq 1$, $\alpha\in (1, 2)$ and $D$ is a $C^{1,1}$
open subset of $\R^n$. Let $\delta_D(x)$ be the Euclidean distance
between $x$ and $D^c$.
\begin{description}
\item{\rm (i)} For every $T>0$, on $(0, T]\times D\times D$
$$
p_D(t, x, y)\asymp t^{-n/\alpha} \left(1\wedge
\frac{t^{1/\alpha}}{|x-y|}\right)^{n+\alpha} \left( 1\wedge
\frac{\delta_D(x)}{t^{1/\alpha}} \right)^{\alpha-1} \left( 1\wedge
\frac{\delta_D(y)}{t^{1/\alpha}}\right)^{\alpha-1}.
$$
\item{\rm (ii)}
Suppose in addition that $D$ is bounded. For every $T>0$, there
exist positive constants $c_1<c_2$ such that for all $(t, x, y)\in
[T, \infty)\times D\times D$,
$$
c_1e^{-\lambda_1 t} \delta_D (x)^{\alpha-1} \delta_D
(y)^{\alpha-1}\leq p_D(t, x, y) \leq c_2 e^{-\lambda_1 t} \delta_D
(x)^{\alpha-1}\delta_D (y)^{\alpha-1},
$$
where $-\lambda_1<0$ is the largest eigenvalue of the $L^2$-generator of $X$.
\end{description}
\end{thm}

\medskip

By integrating the above two-sided heat kernel estimates in Theorem
\ref{T:7.2} with respect to $t$, one can easily obtain the following
sharp two-sided estimate of the Green function
$G_D(x,y)=\int_0^\infty p_D(t, x, y)dt$ of a censored stable process
in a bounded $C^{1,1}$ open set $D$.

\begin{cor}[{\cite[Corollary 1.2]{CKS2}}]\label{C:7.3}
Suppose that $n\geq 1$, $\alpha \in (1, 2)$ and $D$ is a bounded
$C^{1,1}$ open set in $\R^n$.   Then on  $D\times D$, we have
$$
G_D(x, y)\, \asymp\, \begin{cases}  \frac1{|x-y|^{n-\alpha}}\left(
1\wedge  \frac{ \delta_D (x)\delta_D (y)}{|x-y|^2}\right)^{\alpha-1}
\qquad &\hbox{when } n\geq  2,  \\
\big( \delta_D(x)  \delta_D (y)\big)^{(\alpha-1)/2} \wedge
\left(\frac{ \delta_D(x) \delta_D (y)}{ |x-y|} \right)^{\alpha-1}
&\hbox{when } n=1.
\end{cases}
$$
\end{cor}

\bigskip

Sharp two-sided estimates on the Green function are very important
in understanding deep potential theoretic properties of Markov
processes. When $D$ is a bounded $C^{1,1}$ connected open sets in $
\R^n$ and $n\geq 2$, estimates in Corollary \ref{C:7.3} had been
obtained in \cite{CKi}.

Our approach in \cite{CKS2}  is adapted from that of \cite{CKS1}.
In \cite{CKS1}, the following domain monotonicity for the killed
symmetric stable processes is used in a crucial way. Let $Z$ be a
symmetric $\alpha$-stable process and $Z^D$ be the subprocess of $Z$
killed upon leaving an open set $D$. If $U$ is an open subset of
$D$, then $Z^U$ is a subprocess of $Z^D$ killed upon leaving $U$.
However  censored stable-like processes do not have this kind of
domain monotonicity. So there are new challenges to overcome when
studying heat kernel estimates for censored stable processes.
A quantitative version of the intrinsic
ultracontractivity, a crucial use of boundary Harnack inequality
for censored stable process and the reflected stable process $\overline X$
all played an important role in our approach in \cite{CKS2}.

\section{Concluding Remarks}

In this paper, we surveyed some recent progress in the study of fine
potential theoretic properties of various models of symmetric
discontinuous Markov processes that the author is involved. To keep
the exposition as transparent as possible,
sometimes we did not state the results in its most
general form. For example, results in Sections 2 and 3 hold for
general $d$-sets $F$ in $\R^n$ and for $F$ being a measure-metric
space satisfying certain conditions, see \cite{CK1, CK2}; and the
Dirichlet heat kernel estimates in Section 7 in fact holds also
for censored stable-like processes,
see \cite{CKS2}. Two-sided transition density function estimates for
 relativistic stable processes  in $C^{1,1}$ open sets
 have recently been established in \cite{CKS3}.
The study of sharp two-sided heat kernel estimates for discontinuous
Markov processes is in its early stage and is currently a very
active research area. There are many questions waiting to be
answered and several active studies are currently underway. For
instance, it is natural to study the large time estimate for $p(t,
x, y)$ of the processes considered in Section 3 for the case of
 $\gamma_2>\gamma_1>0$.
It is also nature and important to investigate the sharp two-sided
heat kernel estimates for stable processes of mixed type
in $C^{1,1}$-open sets, and for L\'evy processes that
 is the independent sum of a Brownian motion and a symmetric
 stable process in $C^{1,1}$-domains. Some promising progress has
already been made in these studies.

\vskip 0.3truein

\noindent {\bf Zhen-Qing Chen}

\noindent Department of Mathematics, University of Washington,
Seattle, WA 98195, USA \ and \\
Department of Mathematics,  Beijing Institute of Technology,
Beijing 100081, P. R. China

\noindent E-mail: zchen@math.washington.edu

\end{document}